\documentclass[11pt]{article}
\usepackage{cite}
\usepackage{mathrsfs}
\usepackage{amsfonts}
\usepackage{amsmath}
\usepackage{amsfonts,amssymb,color}
\usepackage{dsfont}
\usepackage{curves}
\usepackage{mathrsfs}
\usepackage{pifont}
\usepackage{amssymb}
\usepackage{latexsym,amsmath,amssymb,amsfonts,epsfig,graphicx,cite,psfrag}
\usepackage{eepic,color,colordvi,amscd}

\newtheorem{theorem}{Theorem}[section]

\newtheorem{lemma}{Lemma}[section]
\newtheorem{corollary}{Corollary}[section]

\newtheorem{claim}{Claim}[section]

\newcommand{\qed}{\hfill\rule{0.5em}{0.809em}}

\def\emptyset{\mbox{{\rm \O}}}

\def\qed{\hfill \rule{4pt}{7pt}}

\def\pf{\noindent {\it Proof. }}

\begin{document}
	
	\title{Nearly optimal coloring of some  $C_4$-free graphs}
	\author{Ran Chen\footnote{Email: 1918549795@qq.com},  \; Baogang  Xu\footnote{Email: baogxu@njnu.edu.cn. Supported by NSFC 11931006}\\\\
		\small Institute of Mathematics, School of Mathematical Sciences\\
		\small Nanjing Normal University, 1 Wenyuan Road,  Nanjing, 210023,  China}
	\date{}

\maketitle

\vskip -15pt

\begin{abstract}
A class ${\cal G}$ of graphs is $\chi$-{\em polydet} if ${\cal G}$ has a polynomial binding function $f$ and there is a polynomial time algorithm to determine an $f(\omega(G))$-coloring of $G\in {\cal G}$.
Let $P_t$ and $C_t$ denote a path and a cycle on $t$ vertices, respectively. A {\em bull} consists of a triangle with two disjoint pendant edges, a {\em hammer} is obtained by identifying an end of $P_3$ with a vertex of a triangle, a {\em fork$^+$} is obtained from $K_{1, 3}$ by subdividing an edge twice. Let $H$ be a bull or a hammer, and $F$ be a $P_7$ or a fork$^+$. We determine all $(C_3, C_4, F)$-free graphs without clique cutsets and universal cliques, and present a close relation between $(C_4, F, H)$-free graphs and the Petersen graph. As a consequence, we show that the classes of $(C_4, F, H)$-free graphs are $\chi$-polydet with nearly optimal linear binding functions.
	\begin{flushleft}
	\end{flushleft}
	
\end{abstract}

\section{Introduction}

All graphs considered in this paper are finite and simple. We follow \cite{BM08} for
undefined notations and terminologies. Let $G$ be a graph, let $v\in V(G)$, and let $X$ and $Y$ be two subsets of $V(G)$. Let $E(X, Y)$ be the set of edges with one end in $X$ and the other in $Y$. We say that $v$ is {\em complete} (resp. {\em anticomplete}) to $X$ if $|E(v, X)|=|X|$ (resp. $E(v, X)=\emptyset$), and say that $X$ is complete (resp. anticomplete) to $Y$ if each vertex of $X$ is complete (resp. anticomplete) to $Y$.

Let $N_G(v)$ be the set of vertices adjacent to $v$, $d_G(v)=|N_G(v)|$. Let $N_G(X)=\{u\in V(G)\setminus X\;|\; u$ has a
neighbor in $X\}$, and $N_G[X]=X\cup N_G(X)$. Let $G[X]$ be the subgraph of $G$ induced by $X$. If it does not cause any confusion, we usually omit the subscript $G$. For $u$, $v\in V(G)$, we write $u\sim v$ if $uv\in E(G)$, and write $u\not\sim v$ if $uv\not\in E(G)$. For $X\subset V(G)$ and $v\in V(G)$, let $N_X(v)=N(v)\cap X$.

We say that a graph $G$ contains a graph $H$ if $H$ is isomorphic to an induced subgraph of $G$, and say that $G$ is $H$-{\em free} if it does not contain $H$.
For a family $\{H_1,H_2,\cdots\}$ of graphs, $G$ is $(H_1, H_2,\cdots)$-free if $G$ is
$H$-free for every $H\in \{H_1,H_2,\cdots\}$.

A {\em clique} (resp. {\em stable set}) of $G$ is a set of mutually adjacent
(resp. non-adjacent) vertices in $G$. The {\em clique number} (resp. {\em stability number}) of $G$, denoted by $\omega(G)$ (resp. $\alpha(G)$), is the maximum size of a clique (resp. stable set) in $G$.

Let $H$ be a graph with  $V(H)=\{v_1, v_2, \cdots, v_n\}$. A {\em clique blowup} of $H$ is
any graph $G$ such that $V(G)$ can be partitioned into $n$ cliques,
say $A_1, A_2, \ldots, A_n$, such that $A_i$ is complete to $A_j$ in $G$ if $v_i\sim v_j$ in $H$, and $A_i$ is anticomplete to $A_j$ in $G$ if $v_i\not\sim v_j$ in $H$.
Let $G$ be a clique blowup of $H$. We call $G$ a $t$-{\em clique blowup} of $H$
if $|A_i|=t$ for all $i$, and call $G$ a {\em nonempty clique blowup} of $H$
if $A_i\ne\emptyset$ for all $i$. {\em Under this literature, an induced subgraph of $H$ can be viewed as a clique blowup of $H$ with $A_i$ of sizes 1 or 0}, for all $i$.

Let $k$ be a positive integer. A $k$-{\em coloring } of $G$ is a function $\phi: V(G)\rightarrow \{1,\cdots,k\}$ such that $\phi(u)\ne \phi(v)$ if $u\sim v$. The {\em chromatic number} $\chi(G)$ of $G$ is the minimum number $k$ for which $G$ has a $k$-coloring. A graph is {\em perfect} if all its induced subgraphs $H$ satisfy $\chi(H)=\omega(H)$. An induced cycle of length at least 4 is called  a {\em hole}, and its complement is called an {\em antihole}. A {\em $k$-hole} is a hole of length $k$. A hole or antihole is {\em odd} or {\em even} if it has odd or even number of vertices. In 2006, Chudnovsky {\em et al} \cite{CRST2006} proved the Strong Perfect Graph Theorem. 
\begin{theorem}\label{perfect}{\em\cite{CRST2006}}
	A graph is perfect if and only if it is (odd hole, odd antihole)-free.
\end{theorem}

Let ${\cal F}$ be a hereditary family of graphs.  If there is a function $f$ such that $\chi(H)\leq f(\omega(H))$ for each graph $H$ in ${\cal F}$, then we say that ${\cal F}$ is $\chi$-{\em bounded}, and call $f$ a {\em binding function} of ${\cal F}$ \cite{G75}. Gy\'{a}rf\'{a}s \cite{G75}, and Sumner \cite{T82} independently, conjectured that the class of $T$-free graphs is $\chi$-bounded for any forest $T$.
Interested readers are referred to \cite{RS2004, SR2019, SS18} for its related problems and progresses.

\begin{figure}[htbp]
	\begin{center}
		\includegraphics[width=8cm]{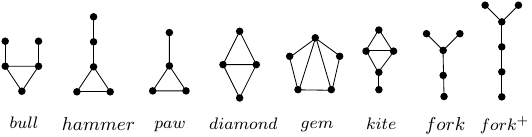}
	\end{center}
	\vskip -25pt
	\caption{Illustration of some special forbidden graphs.}
\label{fig-1}
\end{figure}


Some configurations used in this paper such as {\em paw}, {\em diamond}, {\em gem}, {\em kite, fork}, and {\em fork}$^+$ are shown in Figure~\ref{fig-1}. Since every antihole with at least six vertices has a 4-hole $C_4$, the $\chi$-boundedness of subclasses of $C_4$-free graphs are studies extensively. Choudum {\em et al} \cite{CK2010} proved $\chi(G)\leq\lceil\frac{5}{4}\omega(G)\rceil$ if $G$ is $(P_2\cup P_3, C_4)$-free.  Chudnovsky {\em et al} \cite{CHKK2021} proved  $\chi(G)\leq\frac{3}{2}\omega(G)$ if $G$ is (fork, $C_4$)-free.
Gasper and Huang \cite{GSH17} proved that $\chi(G)\leq \lfloor\frac{3}{2}\omega(G)\rfloor$ for $(P_6, C_4)$-free graph $G$, and Karthick and Maffray \cite{KM19} improved the upper bound to $\lceil\frac{5}{4}\omega(G)\rceil$. Huang \cite{H22} proved that $\chi(G)\leq\lceil\frac{11}{9}\omega(G)\rceil$ if $G$ is $(P_7, C_4, C_5)$-free, which is optimal. Generalize some results of Choudum {\em et al} \cite{CKB21}, the current authors and Wu  proved \cite{CWX2023} that $\chi(G)\leq \max\{3,\omega(G)\}$ if $G$ is ($P_7, C_4$, diamond)-free, and $\chi(G)\leq 2\omega(G)-1$ if $G$ is ($P_7, C_4$, gem)-free, and they also proved that $\chi(G)\leq \omega(G)+1$ if $G$ is ($P_7, C_4$, kite)-free.

Motivated by a new concept {\em Polyanna} introduced by Chudnovsky {\em et al}   \cite{CCDO2023}, and by the above mentioned linear binding function of subclasses of $C_4$-free graphs, we studied $(C_4$, bull)-free graphs and $(C_4$, hammer)-free graphs in \cite{CX2024}, and proved the following conclusion.
\begin{theorem}\label{bull,hammer-1}{\em \cite{CX2024}}
Let ${\cal G}$ be a  $\chi$-bounded class of graphs. Let $H$ be a bull or a hammer, and let $G$ be a connected $(C_4, H)$-free graph which has no clique cutsets or universal cliques. Then,
\begin{itemize}
\item $G$ is a clique blowup of some graph of girth at least $5$ if $H$ is a bull, and $G$ has girth at least 5 if $H$ is a hammer.
\item  $(C_4, H)$-free graphs of ${\cal G}$ are always linearly $\chi$-bounded.
\end{itemize}
\end{theorem}

Note that the paths and stars are the two extremal classes of trees. There are a lot of results on binding functions of $P_t$-free or $K_{1, 3}$-free related graphs.
Gy\'{a}rf\'{a}s \cite{G75} proved that $\chi(G)\leq (t-1)^{\omega(G)-1}$  for all $P_t$-free graphs, and this upper bound was improved to $(t-2)^{\omega(G)-1}$ by Gravier {\em et al} \cite{GHM03}. On $P_5$-free graphs, its best known binding function,due to Scott {\em et al} \cite{SSS2023}, is $\omega(G)^{\log_2(\omega(G))}$ for $\omega(G)\ge3$. Chudnovsky {\em et al} \cite{CS2010} proved that a connected $K_{1,3}$-free graphs $G$ with $\alpha(G)\geq3$ satisfies $\chi(G)\leq2\omega(G)$. Liu {\em et al} \cite{LSWY23} proved that $\chi(G)\le7\omega^2(G)$ if $G$ is fork-free, which answers a problem raised in \cite{KKS2022, RS2004}. Chudnovsky {\em et al} \cite{CHKK2021} proved that $\chi(G)\leq\frac{3}{2}\omega(G)$ if $G$ is (fork, $C_4$)-free. Interested readers can find  more results and problems on binding functions of subclasses of fork-free graphs in \cite{CCS2020,KKS2022, RS2004, WX2023}.

Let $H$ be a bull or a hammer. From Theorem~\ref{bull,hammer-1}, both  $(C_4, P_7, H)$-free graphs and ($C_4$, fork$^+, H)$-free graphs are linearly $\chi$-bounded. We can do better on these two classes of graphs. In this paper, we prove the following Theorems~\ref{fork^+} and \ref{P_7} on the structural decompositions of ($C_4$, fork$^+, H)$-free graphs and $(C_4, P_7, H)$-free graphs, which show that these graphs are closely related to the Petersen graph. Then, we deduce polynomial time algorithms to optimally color these graphs.

\begin{figure}[htbp]
	\begin{center}
		\includegraphics[width=12cm]{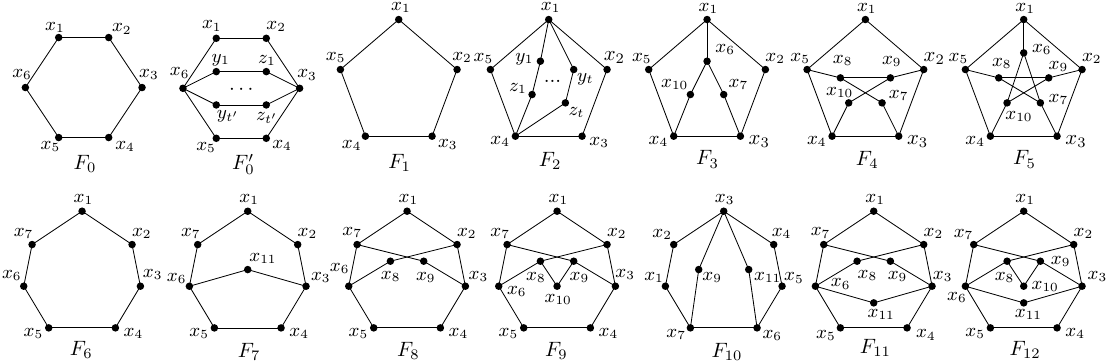}
	\end{center}
	\vskip -25pt
	\caption{Illustration of graphs $F_0\sim F_{12}$.}
\label{fig-2}
\end{figure}

Let $F_0\sim F_{12}$ be the graphs as shown in Figure~\ref{fig-2}, where $F_5$ is just the Petersen graph, and the parameter $t$ and $t'$ appeared in $F_2$ and $F'_{0}$
are positive integers. When $t=1$, we particularly refer to the graph $F_2$ as $F'_2$.
Let $\mathfrak{F} = \{F_0, F'_0, F_1, F_2, \ldots, F_{12}\}$ and, let $\mathfrak{F'} = \{F'_2, F_3, F_4, F_5\}\cup\{C_k~|~k~\mbox{is a positive integer of at least 5}\}$. A {\em clique cutset} of $G$ is a clique $K$ in $G$ such that $G-K$ has more components than $G$, and  a {\em cutvertex} is a clique cutset of size 1. A clique $K$ is call a {\em universal clique} of $G$ if $K$ is complete to $V(G)\setminus K$.

\begin{theorem}\label{fork^+}
Let $H$ be a bull or a hammer, and let $G$ be a connected $(C_4, fork^+, H)$-free graph without clique cutsets or universal cliques. Then, $G\in\mathfrak{F'}$ if $H$ is a hammer,  and $G$ is a nonempty clique blowup of a graph in $\mathfrak{F'}$ if $H$ is a bull.
\end{theorem}
\begin{theorem}\label{P_7}
Let $H$ be a bull or a hammer, and let $G$ be a connected $(C_4, P_7, H)$-free graph without clique cutsets or universal cliques. Then,
$G\in \mathfrak{F}$ if $H$ is a hammer, and
$G$ is a nonempty  clique blowup of a graph in $\mathfrak{F}$ if $H$ is a bull.
\end{theorem}

It is easy to check that $F_5-x_1$ is isomorphic to $F_{12}-\{x_4, x_5\}$. Since $F_5$  is just the Petersen graph, we have that $F_{12}$ can be obtained from the Petersen graph by deleting a vertex and adding a path of length 3 to join some two vertices of degree 3. Since $F_0$, $F_1$, $F'_2$, $F_3$ and $F_4$ are all induced subgraphs of $F_5$, and since $F_6\sim F_{11}$ are all induced subgraphs of $F_{12}$, Theorems~\ref{fork^+} and \ref{P_7} assert that  ($C_4$, fork$^+, H)$-free graphs and $(C_4, P_7, H)$-free graphs are closely related to the Petersen graph, where $H$ is a bull or a hammer.

Notice that all graphs in $\mathfrak{F'}\cup \mathfrak{F'}$ are 3-colorable, and clique cutsets and universal cliques are reducible in coloring of graphs. As an immediate consequence of
Theorems~\ref{fork^+} and \ref{P_7}, one can prove the following Corollary~\ref{hammer-3} by a simple induction on $|V(G)|$.
\begin{corollary}\label{hammer-3}
Let $F$ be a $P_7$ or a fork$^+$, and let $G$ be a $(C_4, F$, hammer)-free graph. Then, $\chi(G)\leq\omega(G)+1$.
\end{corollary}

Let $k$ be a positive integer. We use $G^k$ to denote the $k$-clique blowup of a graph $G$. We can easily verify that both $C_5^k$ and $C_7^k$ are ($C_4, H$, bull)-free for $H$ being a $P_7$ or a fork$^+$, and $C_9^k$ is ($C_4$, fork$^+$, bull)-free. Also, we have that  $\chi(C_5^k)\geq\lceil\frac{|V(C_5^k)|}{\alpha(C_5^k)}\rceil=\lceil\frac{5k}{2}\rceil=\lceil\frac{5}{4}\omega(C_5^k)\rceil$, $\chi(C_7^k)\geq\lceil\frac{|V(C_7^k)|}{\alpha(C_7^k)}\rceil=\lceil\frac{7k}{3}\rceil=\lceil\frac{7}{6}\omega(C_7^k)\rceil$, and $\chi(C_9^k)\geq\lceil\frac{|V(C_9^k)|}{\alpha(C_9^k)}\rceil=\lceil\frac{9k}{4}\rceil=\lceil\frac{9}{8}\omega(C_9^k)\rceil$.
We will show that by excluding $C_5^k$ and $C_7^k$ further,  ($C_4, H$, bull)-free graphs $G$ satisfy $\chi(G)\le \omega(G)+1$ if $H$ is a $P_7$, and $\chi(G)\leq\lceil\frac{9}{8}\omega(G)\rceil$ if $H$ is a fork$^+$.

\begin{theorem}\label{bull-3}
Let $H$ be a $P_7$ or a fork$^+$, and let $G$ be a $(C_4, H$, bull)-free graph. Then,
\begin{itemize}
\item $\chi(G)\leq\lceil\frac{5}{4}\omega(G)\rceil$, and $\chi(G)\leq\lceil\frac{7}{6}\omega(G)\rceil$ if $G$ is further $C_5^2$-free.

\item $\chi(G)\leq\omega(G)+1$ if $H=P_7$ and $G$ is further $(C_5^2, C_7^4)$-free.

\item $\chi(G)\leq\lceil\frac{9}{8}\omega(G)\rceil$ if $H=fork^+$ and $G$ is further $(C_5^2, C_7^4)$-free.
\end{itemize}
Where all the bounds $\lceil\frac{5}{4}\omega(G)\rceil$, $\lceil\frac{7}{6}\omega(G)\rceil$ and $\lceil\frac{9}{8}\omega(G)\rceil$ are reachable.
\end{theorem}

A class ${\cal G}$ of graphs is said to be $\chi$-{\em polydet} \cite{SR2019} if ${\cal G}$ has a polynomial binding function $f$ and there is a polynomial time algorithm to determine an $f(\omega(G))$-coloring for each $G\in {\cal G}$.

Suppose that $G$ has a clique cutset $Q$, and suppose that $V(G)\setminus Q$ is partitioned into two nonempty subsets $V_1$ and $V_2$ that are anticomplete to each other. Let $G_1=G[V_1\cup Q]$ and $G_2=G[V_2\cup Q]$. By repeating this procedure to both $G_1$ and $G_2$ until all the resulted subgraphs have no clique cutsets, $V(G)$ can be partitioned into a collection of subsets of which each induces a subgraph without clique cutsets. We can use a binary tree $T$ to represent this process, and call $T$ a {\em clique-cutsets-decomposition}.

Tarjan \cite{RET1985} showed that for a connected graph with $n$ vertices and $m$ edges, its clique-cutsets-decomposition (not necessarily unique) can be found in $O(nm)$ times. Tarjan also showed that if the maximum weight clique of each subgraph induced by a leaf of $T$ can be found, then one can find the maximum weight clique of $G$ in $O(n^2)$ times.

Chudnovsky {\em et al} \cite{CLMTV2019} showed that an optimal coloring of any  $C_4$-free perfect graph $G$ can be found in $O(|V(G)|^9)$. If $G$ has $t$ maximal cliques, then we can take $O(t|V(G)|^3)$ times to find them \cite{MU2004, TIAS1977}. Since  a $C_4$-free graph $G$ has $O(|V(G)|^2)$ maximal cliques  \cite{A1991, F1989}, we  can find the maximum weight cliques for any $C_4$-free graph $G$ in $O(|V(G)|^5)$ times.  Let $F$ be a $P_7$ or a fork$^+$, and let $H$ be a bull or a hammer. By Theorems~\ref{fork^+} and \ref{P_7}, if $G$ is a $(C_4, F, H)$-free graph, then in its clique-cutsets-decomposition, each leaf represents a subset of vertices which induces a well-defined graph with a universal clique (maybe empty). Applying Tarjan's algorithm,  we can, in at most  $O(|V(G)|^3)$ times, find a nearly optimal coloring and the maximum weight cliques of a $(C_4, F, H$)-free graph. Therefore, we have
\begin{corollary}
Let $F$ be a $P_7$ or a fork$^+$, and let $H$ be a bull or a hammer. Then, $(C_4, F, H$)-free graphs are $\chi$-polydet with linear binding functions.
\end{corollary}

We will prove Theorem~\ref{fork^+} in Section~\ref{3}, prove Theorem~\ref{P_7} in Section~\ref{2}, and prove Theorem~\ref{bull-3} in Section~\ref{4}.

\section{Proof of Theorem~\ref{fork^+}}\label{3}
This section is devoted to prove Theorem~\ref{fork^+}. We begin from (fork$^+,C_4,C_3)$-free graphs.

\begin{lemma}\label{fork^+,C_4,C_3'}
Let $G$ be a connected perfect (fork$^+, C_4, C_3$)-free graph without clique cutsets.
If $G$ is not a complete graph, then $G$ is an even-hole.
\end{lemma}
\pf By Theorem~\ref{perfect}, $G$ is odd-hole-free. If $G$ is even-hole-free, then $G$ is a chordal graph, and so is complete as it has no clique cutsets.
Suppose that $G$ has an even-hole, say  $C=v_1v_2\cdots v_{2q}v_1$, where $q\ge 3$.
If $N(V(C))\ne\emptyset$, let $x\in N(V(C))$ and suppose, without loss of generality, that $x\sim v_1$,
then $x$ is anticomplete to $\{v_2, v_3, v_4, v_{2q}\}$, and so $G[\{x,v_1,v_2,v_3,v_4,v_{2q}\}]$
is a fork$^+$, a contradiction. Therefore, $N(V(C))=\emptyset$, and so
$G$ is an even-hole. \qed

\begin{lemma}\label{fork^+,C_4,C_3}
Let $G$ be a connected imperfect $(fork^+,C_4,C_3)$-free graph without clique cutsets.
If $G$ is not an odd-hole, then $G$ is isomorphic to a graph in \{$F_2', F_3, F_4, F_5$\}.
\end{lemma}
\pf Since $G$ is not perfect, it follows from Theorem~\ref{perfect} that $G$ has an odd-hole. Let $v_1v_2\cdots v_{2q+1}v_1$ be an odd-hole of $G$, where $q\ge 2$, and
let $L=\{v_1,v_2,\cdots,v_{2q+1}\}$ and $R=V(G)\setminus(N(L)\cup L)$.
Since $G$ is not an odd hole, we have that $N(L)\ne\emptyset$. Without loss of generality, suppose that $N(L)$ has a vertex $x$ adjacent to $v_1$.

Suppose that $q\ge 3$. Then, $x$ is anticomplete to $\{v_2, v_3, v_{2q+1}\}$ to avoid a
$C_3$ or $C_4$, and $x\sim v_4$ to avoid an induced fork$^+$
on $\{x,v_1,v_2,v_3,v_4,v_{2q+1}\}$. But then, $x\not\sim v_5$, which forces
$G[\{x,v_1,v_3,v_4,v_5,v_{2q+1}\}]$ to be an induced fork$^+$, a contradiction.
Therefore, $q=2$. During the following proof of Lemma~\ref{fork^+,C_4,C_3}, every subscript is understood to be modulo 5.

Since $G$ has girth at least 5, we have that each vertex of
$N(L)$ has exactly one neighbor in $L$.
Let $i\in\{1,\cdots, 5\}$. We define $X_i=\{x\in N(L)~|~N_L(x)=\{v_i\}\}$,
and let $X=\bigcup_{i=1}^5 X_i$. Then, $N(L)=X$, and $V(G)=L\cup X\cup R$.

If $E(X_i, X_{i+1})\ne\emptyset$, then $x_iv_iv_{i+1}x_{i+1}x_i$ is a 4-hole for some $x_i\in X_i$ and $x_{i+1}\in X_{i+1}$. If $X_i$ is not complete to $X_{i+2}$, then
for some $x'_i\in X_i$ and $x_{i+2}\in
X_{i+2}$, $G[\{x'_i, v_i, v_{i+1}, v_{i+2}, x_{i+2}, v_{i+3}\}]$ is a fork$^+$.
If $X_i$ has two vertices, say  $x$ and $x'$, then $xv_1x'x$ is a triangle if
$x\sim x'$, and $G[\{x,x',v_1,v_2,v_3,v_4\}]$ is a fork$^+$ if $x\not\sim x'$. This shows that
\begin{equation}\label{X}
	\mbox{$|X_i|\leq1$, $X_i$ is anticomplete to $X_{i+1}$ and complete to $X_{i+2}$.}
\end{equation}

In this section, we always assume that $X_i=\{x_i\}$ if $X_i\ne\emptyset$. We show next that
\begin{equation}\label{R'}
	R=\emptyset.
\end{equation}

Suppose $R\ne\emptyset$. We may assume, by symmetry, that
$r\sim x_1$ for some $r\in R$. Let $Q$ be the component of $G[R]$ that contains $r$.
If $X_3\ne\emptyset$, then $x_1\sim x_3$ by (\ref{X}), and so $G[\{r, x_1, x_3, v_1, v_5, v_4\}]$ is a fork$^+$.
If $X_4\ne\emptyset$, then $x_4\sim x_1$ by (\ref{X}), and so $G[\{r, x_1, x_4, v_1, v_2, v_3\}]$ is a fork$^+$.
This shows that $X_3=X_4=\emptyset$. Since $x_1$ is not a cutvertex of $G$, we may assume by symmetry that $X_2=\{x_2\}$ and
$N_{V(Q)}(x_2)\ne\emptyset$. With a similar argument as above used to deal with $X_1$, we can show that $X_{5}=\emptyset$. Then, $\{v_1,v_2\}$ is a clique cutset of $G$. Therefore, (\ref{R'}) holds.

Now, we have that $V(G)=L\cup X$, where $X=N(L)\ne \emptyset$. Suppose, without loss of generality,  that $X_1\ne\emptyset$.
Since $x_1$ is anticomplete to $X_2\cup X_5$ by (\ref{X}),  and since $v_1$ is not a
cutvertex, we may by symmetry assume that $X_3\ne\emptyset$. Then, $x_1\sim x_3$ by
(\ref{X}).

If  $X_4=X_5=\emptyset$, then $X_2=\emptyset$ as otherwise $x_2$ is a cutvertex by (\ref{X}),
and thus $G$ is isomorphic to $F_2'$. Suppose by symmetry that $X_4\ne \emptyset$. By (\ref{X}), $x_1\sim x_4$ and $x_3\not\sim x_4$.
If $X_2\ne\emptyset$ and $X_5\ne\emptyset$, then by (\ref{X}), we have $x_2\sim x_4$,
$x_2\sim x_5$, and $x_3\sim x_5$, and thus  $G$ is isomorphic to $F_5$. If $X_2=X_5=\emptyset$,
then $G$ is isomorphic to $F_3$. If $X_2\ne\emptyset$ and $X_5=\emptyset$, then $x_2\sim x_4$ by (\ref{X}),
and $G$ is isomorphic to $F_4$. The same happens if $X_2=\emptyset$ and $X_5\ne\emptyset$.
Therefore, $G$ is isomorphic to a graph in $\{F_3, F_4, F_5\}$.
This completes the proof of Lemma~\ref{fork^+,C_4,C_3}. \qed

Follows from Lemmas~\ref{fork^+,C_4,C_3'} and \ref{fork^+,C_4,C_3}, we have immediately
the following Lemma~\ref{fork^+,C_4,C_3''}.
\begin{lemma}\label{fork^+,C_4,C_3''}
	Let $G$ be a connected (fork$^+$,$C_4,C_3$)-free graph without clique cutsets. If $G$ is not a complete graph, then $G$ is in $\mathfrak{F'}$.
\end{lemma}

\noindent\textbf{{\em Proof of Theorem~\ref{fork^+}}:}
Let $H$ a bull or a hammer, and	let $G$ be a connected ($C_4$, fork$^+, H)$-free graph without
clique cutsets or universal cliques. By Theorem~\ref{bull,hammer-1}, $G$ is a clique blowup
of some (fork$^+, C_4, C_3)$-free graph if $H$ is a bull, and $G$ is (fork$^+, C_4, C_3)$-free graph if $H$ is a hammer.
If $H$=hammer, then $G$ is in $\mathfrak{F'}$ by Lemma~\ref{fork^+,C_4,C_3''}.

Suppose that $H$ is a bull, and let $G'$ be a (fork$^+,C_4,C_3)$-free graph such that $G$ is a nonempty
clique blowup of $G'$. Then, $G'$ is a non-complete connected graph without clique cutsets,
and so $G'$ is in $\mathfrak{F'}$ by Lemma~\ref{fork^+,C_4,C_3''}. This completes the proof of Theorem~\ref{fork^+}. \qed

\section{Proof of Theorem~\ref{P_7}}\label{2}

The aim of this section is to prove Theorem~\ref{P_7}, which characterizes the structures of
$(P_7, C_4, bull)$-free graphs and $(P_7,C_4, hammer)$-free graphs. First, we discuss  $(P_7,C_4,C_3)$-free graphs.

\begin{lemma}\label{P_7,C_7,C_4,C_3}
Let $G$ be a connected imperfect $(P_7,C_7,C_4,C_3)$-free graph without clique cutsets.
Then $G$ is isomorphic to a graph in $\{F_1, F_2, F_3, F_4, F_5\}$ (see Figure~$\ref{fig-2}$).
\end{lemma}
\pf Since $G$ is not perfect, it follows from Theorem~\ref{perfect} that $G$ contains a 5-hole.  Let $v_1v_2v_3v_4v_5v_1$ be a 5-hole of $G$. Let $L=\{v_1,v_2,\cdots,v_5\}$, and let $R=V(G)\setminus(N(L)\cup L)$. Since $G$ has no triangles and no $C_4$s, we have that each vertex of $N(L)$ has exactly one neighbor in $L$. For each $i\in \{1, 2, \ldots, 5\}$, we define that $A_i=\{x\in N(L)|~N_L(x)=\{v_i\}\}$, and let $A=\bigcup_{i=1}^5 A_i$.
Then, $V(G)=L\cup A\cup R$, and for each $i\in\{1,2,\cdots,5\}$,
\begin{equation}\label{stable}
	\mbox{$A_i$ is a stable set, and $A_i$ is anticomplete to $A_{i+1}\cup A_{i-1}$,}
\end{equation}
where the summation of subindexes is taken modulo 5.

Next, we prove that
\begin{equation}\label{R}
	R=\emptyset.
\end{equation}

Suppose not. Since $G$ is connected, we may assume by symmetry that $a_1\in A_1$ and $r\in R$ such that $a_1\sim r$. Let $Q$ be the component of $G[R]$ which contains $r$.
Since $G$ is $P_7$-free and $ra_1v_1v_2v_3v_4$ is an induced $P_6$, we have that $a_1$ is
complete to $V(Q)$. So, $V(Q)=\{r\}$ because $G$ is triangle-free. If $r$ has a neighbor,
say $a_2$, in $A_2$, then $a_2ra_1v_1v_5v_4v_3$ is an induced $P_7$ by (\ref{stable}). So,
$r$ is anticomplete to $A_2$, and is anticomplete to $A_5$ by symmetry. If $r$ has a neighbor,
say $a_3$, in $A_3$, then $a_1ra_3v_3v_4v_5v_1a_1$ is a 7-hole if $a_1\not\sim a_3$, and
$a_1ra_3a_1$ is a triangle if $a_1\sim a_3$, both are contradictions. Therefore, $r$ is
anticomplete to $A_3$, and is anticomplete to $A_4$ similarly. Since $G$ is $C_4$-free, we have that $a_1$ is the unique neighbor of $r$ in $G$, and so is a cutvertex of $G$, a contradiction. This proves (\ref{R}).

Now, we have that $V(G)=L\cup A$, and $A_i\cup A_{i+1}$ is stable by (\ref{stable}). Let $E_{i, i+2}$ be the set of edges between $A_i$ and $A_{i+2}$. We show next that for each $i\in\{1,2,\cdots,5\}$,
\begin{equation}\label{matching}
	\mbox{$E_{i, i+2}$ is a matching, and if $N_{A_{i+2}}(x)\ne\emptyset$ for some $x\in A_i$ then $N(A_i\setminus\{x\})\cap A_{i-2}=\emptyset$.}
\end{equation}

Let $x\in A_i$. It suffices to prove, by symmetry, that $x$ has at most one neighbor in $A_{i+2}$, and if $x$ has a neighbor in $A_{i+2}$ then $A_i\setminus\{x\}$ has no neighbor in $A_{i-2}$.  If $x$ has two neighbors, say $y$ and $y'$, in $A_{i+2}$, then a 4-hole $xyv_{i+2}y'x$ appears, a contradiction. If $x$ has a neighbor $y$ in $A_{i+2}$, and some vertex $x'\in A_i\setminus\{x\}$ has a neighbor $y'$ in $A_{i-2}$, then a 7-hole $xyv_{i+2}v_{i-2}y'x'v_ix$ appears. Both are contradictions. Therefore, (\ref{matching}) holds.

\begin{equation}\label{A_i''}
\mbox{If $E_{i, i+2}\ne \emptyset$ then $|A_i|=|A_{i+2}|=|E_{i, i+2}|$}.
\end{equation}

Without loss of generality, we take $i=1$, and suppose that $E_{1, 3}\ne \emptyset$. Choose $a_1\sim a_3$ with $a_1\in A_1$ and $a_3\in A_3$. By (\ref{matching}), $N_{A_4}(A_1\setminus\{a_1\})=\emptyset$, and thus $N(A_1\setminus\{a_1\})\subseteq A_3\cup \{v_1\}$. Since $A_1$ is a stable set and $v_1$ is not a cutvertex, we have that each vertex in $A_1\setminus\{a_1\}$ must have a neighbor in $A_3$. Similarly, each vertex in $A_3\setminus\{a_3\}$ must have a neighbor in $A_1$. Therefore, $|A_1|=|A_3|$, and so $E_{1, 3}$ is a matching of size $|A_1|$ by (\ref{matching}). This proves (\ref{A_i''}).

\medskip

Recall that $V(G)=L\cup A$. If $A=\emptyset$ then $G$ is isomorphic to $F_1$. So, we suppose that $A\ne\emptyset$, and suppose, without loss of generality, that $A_1\ne\emptyset$. By (\ref{matching}), $N(A_1)\subseteq A_{3}\cup A_{4}\cup\{v_1\}$.  Since $v_1$ is not a cutvertex, we may assume by symmetry  that $N(A_1)\cap A_3\ne\emptyset$. Then, by (\ref{matching}) and (\ref{A_i''}), $|A_1|=|A_3|=|E_{1, 3}|\geq 1$ and $E_{1, 3}$ is a matching.

Let $a_1\in A_1$ and $a_3\in A_3$ such that $a_1\sim a_3$. We first consider the case that $E_{4, 1}\ne\emptyset$.

\begin{claim}\label{F_3,4,5}
If $E_{4, 1}\ne\emptyset$ then $G$ is isomorphic to a graph in $\{F_3, F_4, F_5\}$.
\end{claim}
\pf Suppose that $A_1$ is not anticomplete to $A_4$. Since $a_1\sim a_3$, we have that
$A_1\setminus\{a_1\}$ is anticomplete to $A_4$ by (\ref{matching}), and $|A_4|=|A_1|=|E_{4, 1}|$
and $|A_1|=|A_3|=|E_{1, 3}|$ by (\ref{A_i''}), with both $E_{1, 3}$ and $E_{4, 1}$ being matchings.
Therefore, $|A_1|=|A_3|=|A_4|=1$. Let $A_4=\{a_4\}$. Then,
$V(G)=\{a_1, a_3, a_4\}\cup A_2\cup A_5\cup L$ such that $G[\{a_1, a_3, a_4\}]=a_3a_1a_4$.

First, we suppose that $E_{3, 5}\ne\emptyset$. By (\ref{matching}) and
(\ref{A_i''}), $|A_5|=|A_3|=|E_{3, 5}|=1$. Let $A_5=\{a_5\}$. Now, we have that
$V(G)=\{a_1,a_3,a_4,a_5\}\cup A_2\cup L$ such that $G[\{a_1, a_3, a_4, a_5\}]=a_5a_3a_1a_4$
by (\ref{stable}). If $A_2=\emptyset$, then $G$ is isomorphic to $F_4$. Otherwise, let $a_2\in A_2$. By (\ref{stable}),
$N(a_2)\subseteq \{a_4, a_5, v_2\}$, and so  $a_2$ must have a neighbor in $\{a_4, a_5\}$
since $v_2$ cannot be a cutvertex. Without loss of generality, we suppose
$a_2\sim a_4$. Then, $a_2\sim a_5$ for avoiding an induced  $P_7=v_2a_2a_4a_1a_3a_5v_5$,
and thus $|A_2|=|A_5|=|E_{5, 2}|=1$ by (\ref{A_i''}), which forces that $G$ is isomorphic to $F_5$.

Similarly, we may deduce that $G$ is isomorphic to $F_4$ or $F_5$ if $E_{2, 4}\ne\emptyset$. Next, we suppose that $E_{2, 4}=E_{3, 5}=\emptyset$.
\medskip

If $A_2\ne\emptyset$, let $a_2\in A_2$, then $N(a_2)\subseteq \{v_2\}\cup A_5$ by (\ref{stable}), and so $a_2$ must have a neighbor, say $a_5$, in $A_5$ as otherwise $v_2$ is a cutvertex.  But now, $v_2a_2a_5v_5v_4a_4a_1$ is an induced $P_7$ as
$a_5\not\sim a_1$ and $a_5\not\sim a_4$ by (\ref{stable}). A similar contradiction happens
if $A_5\ne \emptyset$. Therefore, $A_2=A_5=\emptyset$, and $V(G)=\{a_1,a_3,a_4\}\cup L$.
Now we have that $G$ is isomorphic to $F_3$. This proves Claim~\ref{F_3,4,5}. \qed

\medskip

By Claim~\ref{F_3,4,5} and by symmetry, we may suppose that $E_{4, 1}=E_{3, 5}=\emptyset$. If $A_4\ne\emptyset$, let $a_4\in A_4$, then by (\ref{stable}), and by the fact that $v_4$ cannot be a cutvertex,  $a_4$ must have a neighbor, say $a_2$,  in $A_2$ which forces an induced $P_7=v_4a_4a_2v_2v_1a_1a_3$,  a contradiction. So, $A_4=\emptyset$. Similarly, $A_5=\emptyset$, and thus $A_2=\emptyset$ as otherwise $v_2$ is a cutvertex. Now, we have that $V(G)=A_1\cup A_3\cup L$, $|A_1|=|A_3|\geq1$, and $E_{1, 3}$ is a matching of size $|A_1|$. Therefore, $G$ is isomorphic to $F_2$. This completes the proof of Lemma~\ref{P_7,C_7,C_4,C_3}. \qed

\begin{lemma}\label{P_7,C_4,C_3}
	Let $G$ be a connected $(P_7,C_4,C_3)$-free graph without clique cutsets.
If $G$ contains a $7$-hole, then $G$ is isomorphic to a graph in \{$F_6, F_7,\cdots, F_{12}$\}.
\end{lemma}
\pf Let $v_1v_2\cdots v_7v_1$ be a 7-hole of $G$.
Let $L=\{v_1,v_2,\cdots,v_7\}$ and $R=V(G)\setminus (N(L)\cup L)$.
For each $i\in\{1,2,\cdots,7\}$, we define $X_i=\{x\in N(L)~|~N_L(x)=\{v_i,v_{i+3}\}\}$.
During the proof of Lemma~\ref{P_7,C_4,C_3}, every subscript is understood to be modulo 7.
Let $X=\bigcup_{i=1}^7 X_i$.

Let $x\in N(L)$. Without loss of generality, suppose that $x\sim v_1$.
Since $G$ has girth at least 5, we have that $x$ is anticomplete to $\{v_2,v_3,v_6,v_7\}$, and thus
must have a neighbor in $\{v_4,v_5\}$ to avoid an induced $P_7=xv_1v_2v_3v_4v_5v_6$.
Obviously, $x$ has exactly one neighbor in $\{v_4,v_5\}$, and hence $x\in X_1$ if $x\sim v_4$
and $x\in X_5$ if $x\sim v_5$. Therefore, $N(L)=X$, and $V(G)=L\cup X\cup R$.

Since $G$ has girth at least 5, we have that $|X_i|\leq1$ for each $i\in\{1,2\cdots,7\}$.
In this section, we always assume that
\begin{equation}\label{X_i}
	\mbox{$X_i=\{x_i\}$ if $X_i\ne\emptyset$.}
\end{equation}

Let $i\in \{1,2,\cdots,7\}$.

If $x_i\sim y$ for some $y\in X_{i+1}\cup X_{i+2}$, then $G$ has a
4-hole $x_iyv_{i+1}v_ix_i$ if $y=x_{i+1}$, and has a triangle $x_iyv_{i+3}x_i$ if $y=x_{i+3}$. If $X_i\ne\emptyset$ and $X_{i+2}\ne\emptyset$, then $G$ has a 4-hole $x_ix_{i+2}v_{i+2}v_{i+3}x_i$ if $x_i\sim x_{i+2}$, and has an induced $P_7=x_iv_iv_{i+1}v_{i+2}x_{i+2}v_{i+5}v_{i+4}$ otherwise. Both are contradictions. Therefore,
\begin{equation}\label{X_i'}
	\mbox{$E(X_i, X_{i+1}\cup X_{i+3})=\emptyset$, and either $X_i=\emptyset$ or $X_{i+2}=\emptyset$.}
\end{equation}

Next, we show that
\begin{equation}\label{X_i'''}
	\mbox{if $X_i\ne\emptyset$ and $X_{i+1}\ne\emptyset$, then $N_{R}(x_i)=N_{R}(x_{i+1})$
and $|N_R(x_i)|\leq1$.}
\end{equation}

Without loss of generality, we set $i=1$. Then,  $X_1=\{x_1\}$ and $X_2=\{x_2\}$ by (\ref{X_i}), and $x_1\not\sim x_{2}$ by (\ref{X_i'}).
Let $r\in N_R(x_1)$. Then,
$x_{2}\sim r$ to avoid an induced $P_7=rx_1v_1v_{7}v_{6}v_{5}x_{2}$.
This implies that $N_R(x_1)\subseteq N_R(x_{2})$. Similarly, we can show $N_R(x_{2})\subseteq N_R(x_1)$,
and thus $N_{R}(x_1)=N_{R}(x_{2})$. Since $G$ is triangle-free, $N_R(x_1)$ must be a stable
set, and any two distinct vertices $r$ and $r'$ in $N_R(x_1)$ would produce a 4-hole $rx_1r'x_{2}r$.
Therefore, $|N_R(x_1)|\leq1$. This proves (\ref{X_i'''}).

Recall that $V(G)=L\cup X\cup R$. We may assume that $X\ne\emptyset$ as otherwise $G$ is  isomorphic to $F_6$. Without loss of generality, suppose $X_1=\{x_1\}$.
If $X=\{x_1\}$, then $R=\emptyset$ as otherwise $x_1$ is a cutvertex, and so $G$ is isomorphic to $F_7$.
So, we further suppose that $X\ne \{x_1\}$.

Since $X_1\ne \emptyset$, by (\ref{X_i'}), we have that $X_3=X_6=\emptyset$, and
$\{X_2, X_4, X_5, X_7\}$ contains at most two nonempty sets. By symmetry, we may assume that $X_4=X_7=\emptyset$. Then, $X$ equals one of the set
in $\{\{x_1, x_2\}, \{x_1, x_5\}, \{x_1, x_2, x_5\}\}$. We deal with the three cases separately.

\begin{claim}\label{F_8,9}
 If $X=\{x_1,x_2\}$, then $G$ is isomorphic to $F_8$ or $F_{9}$.
\end{claim}
\pf Suppose that $X=\{x_1,x_2\}$. Then, $V(G)=\{x_1,x_2\}\cup R\cup L$, and $x_1\not\sim x_2$ by (\ref{X_i'}).
If $R=\emptyset$, then $G$ is isomorphic to $F_8$. So, we suppose that $R\ne\emptyset$.
Since $G$ is connected, we may suppose that $N_R(x_1)=N_R(x_2)=\{r\}$ by (\ref{X_i'''}).
If $R\ne \{r\}$, we may choose an $r'\in N_R(r)$, then $G$ has an induced $P_7=r'rx_1v_1v_7v_6v_5$.
Therefore, $R=\{r\}$, and $G$ is isomorphic to $F_{9}$. \qed

\begin{claim}\label{F_{10}}
	If $X=\{x_1,x_5\}$, then $G$ is isomorphic to $F_{10}$.
\end{claim}
\pf Suppose that $X=\{x_1, x_5\}$. Then, $V(G)=\{x_1,x_5\}\cup R\cup L$, and
$x_1\not\sim x_5$ by (\ref{X_i'}). We show that $R=\emptyset$. Suppose not,
we may  assume by symmetry that
$N_R(x_1)\ne\emptyset$. Let $r\in N_R(x_1)$, and let $Q$ be the component of $G[R]$
which contains $r$. To avoid an induced $P_7$ which contains $rx_1v_1v_7v_6v_5$,
we have that  $x_1$ must be complete to $V(Q)$, and hence $V(Q)=\{r\}$ since $G$ is triangle-free.
But then, $r\not\sim x_5$ to avoid a 4-hole $rx_1v_1x_5r$, which forces $x_1$ to be a
cutvertex $G$, contradicting the choice of $G$. Therefore, $R=\emptyset$, and
$G$ must be isomorphic to $F_{10}$. \qed

\begin{claim}\label{F_{11,12}}
	If $X=\{x_1, x_2, x_5\}$, then $G$ is isomorphic to $F_{11}$ or $F_{12}$.
\end{claim}
\pf Suppose that $X=\{x_1, x_2, x_5\}$. Then, $V(G)=\{x_1,x_2,x_5\}\cup R\cup L$,
and $\{x_1,x_2,x_5\}$ is a stable set by (\ref{X_i'}). If $R=\emptyset$ then $G$ is
isomorphic to $F_{11}$. So, we suppose that $R\ne\emptyset$.

We show first that $N_R(x_5)=\emptyset$. Suppose to its contrary,
let $r\in N_R(x_5)$, and let $Q$ be the component of $G[R]$ which contains $r$. Then,
$x_5$ must be complete to $V(Q)$ to avoid an induced $P_7$ which contains $rx_5v_5v_4v_3v_2$,
and hence $V(Q)=\{r\}$ since $G$ is triangle-free.
But now, $r\not\sim x_1$ to avoid a 4-hole $rx_1v_1x_5r$, and $r\not\sim x_2$ to avoid a 4-hole
$rx_2v_5x_5r$, which imply that $x_5$ is a cutvertex of $G$, contradicting the choice of $G$.
Therefore,  $N_R(X_5)=\emptyset$, and so $N_R(x_1)\cup N_R(x_2)\ne\emptyset$.

By (\ref{X_i'''}), we suppose that $N_R(x_1)=N_R(x_2)=\{r'\}$.
If $R\ne \{r'\}$, we may choose $r''\in N_R(r')$,
then $G$ has an induced $P_7=r''r'x_1v_1v_7v_6v_5$, a contradiction.
Therefore, $R=\{r'\}$, and $G$ is isomorphic to $F_{12}$. \qed

Lemma~\ref{P_7,C_4,C_3} follows immediately from Claims~\ref{F_8,9},
\ref{F_{10}} and \ref{F_{11,12}}. \qed

\begin{lemma}\label{P_7,C_4,C_3'}
Let $G$ be a connected perfect ($P_7,C_4,C_3$)-free graph without clique cutsets. If $G$ is not a complete graph, then $G$ is isomorphic to either $F_{0}$ or $F'_{0}$.
\end{lemma}
\pf Since $G$ is perfect, we have that $G$ is ($C_5,C_7$)-free by Theorem~\ref{perfect}.
If $G$ is $C_6$-free, then $G$ is a chordal graph, which must be complete of order  1 or 2
as it is triangle-free and  has no clique cutsets. So, we suppose that $G$ has a 6-hole, say $v_1v_2\cdots v_6v_1$.
Let $L=\{v_1, \cdots, v_6\}$ and $R=V(G)\setminus(N(L)\cup L)$.  During the proof of Lemma~\ref{P_7,C_4,C_3'}, every subscript is
understood to be modulo 6. Since $G$ has girth at least 6, we have  that each vertex of $N(L)$ has exactly one neighbor in $L$. For each
$i\in\{1,2,\cdots,6\}$, we define $Y_i=\{u\in N(L)~|~N(u)\cap L=\{v_i\}\}$, and let
$Y=\bigcup_{i=1}^6 Y_i$. Then, $V(G)=L\cup Y\cup R$.

If $R\ne \emptyset$, let $r\in R$ and suppose by symmetry that $r\sim y_1$ for some $y_1\in Y_1$, then $ry_1v_1v_2v_3v_4v_5$ is an induced $P_7$.
Therefore, $R=\emptyset$, and hence $V(G)=L\cup Y$.

If $Y=\emptyset$ then $G$ is isomorphic to $F_{0}$. Suppose so that
$Y\ne\emptyset$, and suppose, without loss of generality, that $Y_1\ne\emptyset$.

Let $i\in\{1,2,\cdots,6\}$. If $Y_i\ne \emptyset$ and $Y_{i+2}\ne\emptyset$, let
$y_i\in Y_i$ and $y_{i+2}\in Y_{i+2}$, then $y_iv_iv_{i+1}v_{i+2}y_{i+2}y_i$ is an
induced $C_5$ if $y_i\sim y_{i+2}$, and $y_iv_iv_{i+5}v_{i+4}v_{i+3}v_{i+2}y_{i+2}$
is an induced $P_7$ if $y_i\not\sim y_{i+2}$. Therefore,
either $Y_i=\emptyset$ or $Y_{i+2}=\emptyset$.

Since $Y_1\ne \emptyset$ by our assumption, we have $Y_3\cup Y_5=\emptyset$. If $Y_2\ne\emptyset$,
then $Y_4\cup Y_6=\emptyset$, which implies that $\{v_1,v_2\}$ is a clique cutsets.
So, $Y_2=\emptyset$, and $Y_6=\emptyset$ by symmetry. Now, we have that
$Y=Y_1\cup Y_4$. Since both $v_1$ and $v_4$ cannot be cutvertices, and since $Y_1$ and $Y_4$
are both stable, we have that each vertex in $Y_1$ must have a neighbor in $Y_4$,
and each vertex in $Y_4$ must have a neighbor in $Y_1$.
If some vertex $y_1$ of $Y_1$ has two neighbors, say $y$ and $y'$, in $Y_4$, then
$y_1yv_4y'y_1$ would be a 4-hole. A similar contradiction happens if some vertices
of $Y_4$ has two neighbors in $Y_1$. Therefore, $E(Y_1, Y_4)$ is a matching of size $|Y_1|$, and $|Y_1|=|Y_4|$. Now, $G$ is isomorphic to $F'_{0}$.
This completes the proof of Lemma~\ref{P_7,C_4,C_3'}. \qed

Combine Lemmas~\ref{P_7,C_7,C_4,C_3}, \ref{P_7,C_4,C_3} and
\ref{P_7,C_4,C_3'}, we have the following Lemma~\ref{P_7,C_4,C_3''} immediately.
\begin{lemma}\label{P_7,C_4,C_3''}
	Let $G$ be a connected $(P_7,C_4,C_3)$-free graph without clique cutsets.
If $G$ is not a complete graph, then $G$ is in $\mathfrak{F}$.
\end{lemma}

\noindent\textbf{{\em Proof of Theorem~\ref{P_7}}:}  Let $H$ be a bull or a hammer, and
let $G$ be a connected $(P_7,C_4,H)$-free graph without clique cutsets or universal cliques.
By Theorem~\ref{bull,hammer-1}, $G$ is a clique blowup of some $(P_7,C_4,C_3)$-free graph
if $H$ is a bull, and $G$ is $(P_7,C_4,C_3)$-free if $H$ is a hammer.
If $H$ is a hammer, then $G$ is in $\mathfrak{F}$ by Lemma~\ref{P_7,C_4,C_3''}.

Suppose that $H$ is a bull, and suppose that $G$ is a nonempty clique blowup of a graph $G'$.
Since $G$ is connected and has no clique cutsets or universal cliques, we have that
$G'$ is connected, not complete, and has no clique cutsets.
Therefore, $G'$ is in $\mathfrak{F}$ by Lemma~\ref{P_7,C_4,C_3''}.

This completes the proof of Theorem~\ref{P_7}. \qed

\section{Chromatic bound for clique blowups}\label{4}

In this section, we will prove Theorem~\ref{bull-3}. Before that, we first color
some special clique blowups. Throughout this section, we will use the following notations.
Suppose that $G$ is a clique  blowup of $H$ with $V(H)=\{x_1,x_2,\cdots,x_{|V(H)|}\}$.
Let $V(G)$ be partitioned into $|V(H)|$ cliques
$Q_{x_1}$, $Q_{x_2}$, $\cdots$, $Q_{x_{|V(H)|}}$, such that $Q_{x_i}$ is the clique of $G$
corresponding to the vertex  $x_i$ of $H$. We simply write $Q_{x_{i_1}, x_{i_2}, \cdots, x_{i_j}}$ instead of $Q_{x_{i_1}}\cup Q_{x_{i_2}}\cup\cdots\cup Q_{x_{i_j}}$ for $\{i_1, i_2, \cdots, i_j\}\subseteq\{1,2,\cdots,|V(H)|\}$.
We say that two sets $A$ and $B$ {\em meet} if $A\cap B\ne\emptyset$.
Let $G$ be a graph and $S\subseteq V(G)$ be a stable set, we say $S$ is {\em good}
if it meets every clique of size $\omega(G)$.

The following results (Lemmas~\ref{KM}$\sim$ \ref{2-blowup'}) will be used frequently in the proof of Theorem~\ref{bull-3}.
\begin{lemma}\label{KM}{\em\cite{KM19}}
Let $G$ be a graph such that every proper induced subgraph $G'$ of $G$ satisfies
$\chi(G')\leq\lceil\frac{5}{4}\omega(G')\rceil$. Suppose that $G$ has a vertex of degree at most $\lceil\frac{5}{4}\omega(G)\rceil-1$, or has a good stable set, or has a stable set $S$ such that $G-S$ is perfect.
Then, $\chi(G)\leq\lceil\frac{5}{4}\omega(G)\rceil$.
\end{lemma}

In \cite{KM19}, Karthick and Maffray proved that $\chi(G)\leq \lceil \frac{5}{4}\omega(G) \rceil$ if $G$ is a clique blowup of the Petersen graph. Since $F_1$, $F'_2$, $F_3$, and $F_4$ are all induced subgraphs of the Petersen graph $F_5$, we have the following lemma directly.
\begin{lemma}\label{Petersen}
$\chi(G)\leq \lceil \frac{5}{4}\omega(G) \rceil$ if $G$ is a clique blowup of some graph in $\{F_1, F'_2, F_3, F_4, F_5\}$.
\end{lemma}

\begin{lemma}\label{perfect'}{\em\cite{L1972}}
Clqiue blowups of perfect graphs are perfect.
\end{lemma}

Recall that we use $G^k$ to denote the nonempty clique blowup of $G$ such that each vertex is blewup into a clique of size $k$.

\begin{lemma}\label{2-blowup}
	$\chi(F_2^2)=\chi(F_{12}^2)=5$.
\end{lemma}
\pf It is certain that $\chi(F_2^2)\geq\lceil\frac{|V(F_2^2)|}{\alpha(F_2^2)}\rceil
=\lceil\frac{4t+10}{t+2}\rceil=5$, and
$\chi(F_{12}^2)\geq\lceil\frac{|V(F_{12}^2)|}{\alpha(F_{12}^2)}\rceil=5$.

In $F_2^2$, we may color $Q_{x_i}$, $i\in\{1,2,3,4,5\}$, with $\{2i-1, 2i\}$ (taken modulo 5), and color $Q_{y_i}$ and $Q_{z_i}$, $i\in\{1,2,\cdots, t\}$,
with $\{3, 4\}$ and $\{1, 5\}$, respectively. Hence, $\chi(F_2^2)=5$.

For $F_{12}^2$, we can define a 5-coloring $f$ as follows: $f(Q_{x_1})=\{1, 2\}$,  $f(Q_{x_2,x_4})=\{3, 4\}$, $f(Q_{x_3})=\{1, 5\}$, $f(Q_{x_5,x_8})=\{2, 5\}$, $f(Q_{x_6})=\{1, 4\}$, $f(Q_{x_7})=\{3, 5\}$, $f(Q_{x_9})=\{2, 4\}$, $f(Q_{x_{10}})=\{1, 3\}$, $f(Q_{x_{11}})=\{2, 3\}$. This proves Lemma~\ref{2-blowup}. \qed

\begin{lemma}\label{delta}
Let $H$ be a graph, $x\in V(H)$ with $N_H(x)=\{y, z\}$. Let $G$ be a clique blowup of $H$ with $\delta(G)\geq\omega(G)+1$.
Then, $|Q_x|\leq\omega(G)-2$, and $\min\{|Q_y|, |Q_z|\}\geq 2$ unless $Q_x=\emptyset$.
\end{lemma}
\pf It is certain that $\omega(G)\geq2$ since $\delta(G)\geq\omega(G)+1$. If $Q_x=\emptyset$, then we are done. Suppose $Q_x\ne\emptyset$. Then, $\omega(G)\ge \max\{|Q_x|+|Q_y|, |Q_x|+|Q_z|\}\ge |Q_x|+1$.
If $|Q_x|\geq\omega(G)-1$, then $|Q_x|=\omega(G)-1$ and $\max\{|Q_y|,|Q_z|\}=1$, and so
$d_G(x)=|Q_x|+|Q_y|+|Q_z|-1\leq\omega(G)$, a contradiction.
Therefore, $|Q_x|\leq\omega(G)-2$.
Since $\omega(G)\ge \max\{|Q_x|+|Q_y|, |Q_x|+|Q_z|\}$, we have that $\omega(G)+1\leq d_G(x)=|Q_x|+|Q_y|+|Q_z|-1
\leq \min\{\omega(G)-1+|Q_y|,\omega(G)-1+|Q_z|\}$, which implies that $\min\{|Q_{y}|,|Q_z|\}\geq2$. \qed

\begin{lemma}\label{2-blowup'}
Let $H$ be a triangle-free graph with $\chi(H^2)\leq5$,
and let $G$ be a clique blowup of $H$ such that $\chi(G')\leq\lceil\frac{5}{4}\omega(G')\rceil$
for every proper induced subgraph $G'$ of $G$.
Suppose that $\min\{|Q_x|, |Q_y|\}\geq2$ for every edge $xy\in E(H)$ with
$|Q_{x,y}|=\omega(G)$. Then, $\chi(G)\leq\lceil\frac{5}{4}\omega(G)\rceil$.
\end{lemma}
\pf We may suppose that $G$ is a connected imperfect graph. So, $\omega(G)\ge \omega(H)=2$. If $\omega(G)=2$, then $G$ is an induced subgraph of $H$ and $|Q_{x, y}|=2=\omega(G)$ for every
edge $xy\in E(G)\subseteq E(H)$. If $\omega(G)=3$, then $|Q_{x, y}|=3=\omega(G)$ for some edge $xy\in E(H)$  since $H$ is triangle-free. Both contradict $\min\{|Q_x|, |Q_y|\}\geq2$ while $|Q_{x, y}|=\omega(G)$. Therefore, $\omega(G)\geq4$.

Let $K$ be a maximal clique of $G$. Then, $K=Q_{x, y}$ for some edge $xy\in E(H)$ since $H$ is triangle-free. Choose a subset $T$ of $V(G)$ such that $|T\cap Q_v|=\min\{2, |Q_v|\}$ for each $v\in V(H)$. If $|K|=\omega(G)$, then $|K\setminus T|=|K|-4=\omega(G)-4$ because $\min\{|Q_x|,|Q_y|\}\geq2$. Suppose that $|K|\leq\omega(G)-1$. If $\max\{|Q_x|, |Q_y|\}\ge 2$, then $|T\cap Q_{x, y}|\ge 3$, and thus $|K\setminus T|\leq\omega(G)-4$. If $|Q_x|=|Q_y|=1$, then $|K\setminus T|=0\leq\omega(G)-4$. Therefore, 	$\omega(G-T)\leq\omega(G)-4$.

Since $\chi(G[T])\leq \chi(H^2)\leq 5$, we have that
$\chi(G)\leq\chi(G-T)+\chi(G[T])\leq\lceil\frac{5}{4}\omega(G-T)\rceil+5
\leq\lceil\frac{5}{4}(\omega(G)-4)\rceil+5= \lceil\frac{5}{4}\omega(G)\rceil$.
This proves Lemma~\ref{2-blowup'}. \qed

\medskip

Before proving Theorem~\ref{bull-3}, we present a number of lemmas of which all are proved by induction on the number of vertices. In the following proofs, we always suppose that
\begin{equation}\label{imperfect000}
\mbox{$G$ is connected, imperfect and has no clique cutsets}.
\end{equation}

We begin from the clique blowups of odd holes.

\begin{lemma}\label{odd-hole'}
Let $q\ge 2$, let $C=v_1v_2\cdots v_{2q+1}v_1$, and let $G$ be a clique blowup of $C$.
Then, $\chi(C^k)=\lceil\frac{2q+1}{2q}\omega(C^k)\rceil$ for every positive integer $k$,
and $\chi(G)\leq\lceil\frac{2q+1}{2q}\omega(G)\rceil$.
\end{lemma}
\pf Let $k$ be a positive integer. Then, $\omega(C^k)=2k$, and
$\chi(C^k)\geq\lceil\frac{|V(C^k)|}{\alpha(C^k)}\rceil=\lceil\frac{2q+1}{q}k\rceil=
\lceil\frac{2q+1}{2q}\omega(C^k)\rceil$. We can construct a $\lceil\frac{2q+1}{q}k\rceil$-coloring of $C^k$ by coloring $Q_{v_i}$ with $\{(i-1)k+1, (i-1)k+2, \cdots, ik\}$, where the colors are taken modulo
$\lceil\frac{2q+1}{q}k\rceil$. Therefore, $\chi(C^k)=\lceil\frac{2q+1}{q}k\rceil$.

Next, we prove $\chi(G)\leq\lceil\frac{2q+1}{2q}\omega(G)\rceil$ by induction on
$|V(G)|$. If there exists $Q_{v_i}=\emptyset$ for some $i$,
then $G$ is a clique blowup of a path, and the lemma holds by
Lemma~\ref{perfect'}. So, suppose that $Q_{v_i}\ne\emptyset$
for all $i$. If $Q_{v_1, v_2}$ is not a maximum clique of
$G$, let $S=\{v_3,v_5,\cdots,v_{2q+1}\}$, then $S$  is a good stable set of $G$,
and so $\chi(G)\leq\chi(G-S)+1\leq\lceil\frac{2q+1}{2q}(\omega(G)-1)\rceil+1\leq
\lceil\frac{2q+1}{2q}\omega(G)\rceil$ by induction.
So, we may suppose by symmetry that $Q_{v_i, v_{i+1}}$ is a maximum clique of $G$ for all
$i$. Now, $G$ is a $|Q_{v_1}|$-clique blowup
of $C$, and thus $\chi(G)=\lceil\frac{2q+1}{2q}\omega(G)\rceil$.
This proves Lemma~\ref{odd-hole'}. \qed

\begin{lemma}\label{F_2}
Let $G$ be a clique blowup of $F_2$ (see Figure~\ref{fig-2}). Then, $\chi(G)\leq\lceil\frac{5}{4}\omega(G)\rceil$.
\end{lemma}
\pf The lemma holds trivially if $G$ is an induced subgraph of $F_2$. By (\ref{imperfect000}), we have  $Q_{x_i}\ne\emptyset$ if $i\in\{1,4,5\}$.
If $|Q_{x_i}|=1$ for some $i\in\{1,4,5\}$, then $G-\{x_i\}$ is perfect,
and the lemma follows from Lemma~\ref{KM}. So, suppose that
$|Q_{x_i}|\geq2$ if $i\in\{1,4,5\}$.

To avoid a good stable set $S\subseteq\{y_1, y_2, \cdots, y_t, x_2, x_4\}$,
$Q_{x_1, x_5}$ must be a maximum clique of $G$. Similarly, $Q_{x_4, x_5}$ is a
maximum clique of $G$ also. Let $|Q_{x_5}|=a$ and $|Q_{x_1}|=|Q_{x_4}|=\omega(G)-a$. Then,
$a\geq2$ and $\omega(G)-a\geq2$.
We can deduce similarly that, for each $i\in\{1,2,\cdots,t\}$,
$|Q_{y_i}|\geq2$ if $|Q_{x_1, y_i}|=\omega(G)$,
and $|Q_{z_i}|\geq2$ if $|Q_{x_4, z_i}|=\omega(G)$.
So, if $|Q_u|=1$ or $|Q_v|=1$ for some maximum clique $Q_{v, u}$,
then we may by symmetry suppose that $v=y_1$ and $u=z_1$.

Suppose $|Q_{y_1}|=\omega(G)-1$. Then,
$\omega(G)\ge |Q_{x_1, y_1}|=\omega(G)-a+\omega(G)-1=2\omega(G)-a-1$,
and thus $\omega(G)-a=1$ because $\omega(G)>a$, which contradicts
that $|Q_{x_4}|=\omega(G)-a\geq2$. A similar contradiction happens if $|Q_{z_1}|=\omega(G)-1$.
Therefore, for each $uv\in E(F_2)$,
if $|Q_{v, u}|=\omega(G)$ then $\min\{|Q_u|, |Q_{v}|\}\geq 2$.
By Lemmas~\ref{2-blowup} and \ref{2-blowup'}, we have that
$\chi(G)\leq\lceil\frac{5}{4}\omega(G)\rceil$. \qed

\medskip

By Lemma~\ref{Petersen}, we have that $\chi(G)\leq\lceil\frac{5}{4}\omega(G)\rceil$ if $G$ is a clique blowup of some graph in $\{F_1, F_2', F_3, F_4, F_5\}$. By Lemma~\ref{odd-hole'}, we have that $\chi(G)\leq\lceil\frac{7}{6}\omega(G)\rceil<\lceil\frac{5}{4}\omega(G)\rceil$ if $G$ is a clique blowup of $F_6$. Our next lemma deals with the cases for $F_7$ and $F_8$.

\begin{lemma}\label{F_8}
Let $G$ be a clique blowup of $F_8$. Then $\chi(G)\leq\lceil\frac{5}{4}\omega(G)\rceil$.
\end{lemma}
\pf We take $F_8$ to show the proving procedure. It holds trivially if $G$ is an induced subgraph of $F_8$.
We may suppose $\delta(G)\geq\lceil\frac{5}{4}\omega(G)\rceil$ by Lemma~\ref{KM}.
By (\ref{imperfect000}), $Q_{x_1}, Q_{x_2}$, and $Q_{x_7}$ are all nonempty.
By Lemma~\ref{KM}, we may suppose that for each $i\in\{1,2,7\}$,
$|Q_{x_i}|\geq2$ as otherwise $G-\{x_i\}$ is perfect. Then, $|Q_{x_j}|\leq\omega(G)-2$ for $j\in\{1,2,3,6,7,8,9\}$.
Since $\delta(G)\geq\lceil\frac{5}{4}\omega(G)\rceil\geq\omega(G)+1$,
we have that $\max\{|Q_{x_4}|, |Q_{x_5}|\}\leq\omega(G)-2$ by Lemma~\ref{delta}.

Now, we can easily deduce that, for each edge  $uv\in E(F_8)$, $\min\{|Q_v|, |Q_u|\}\geq2$
if $|Q_{u, v}|=\omega(G)$. Lemma~\ref{F_8} follows immediately from Lemmas~\ref{2-blowup}
and \ref{2-blowup'}. \qed

\begin{lemma}\label{F_9}
Let $G$ be a clique blowup of $F_9$. Then $\chi(G)\leq\lceil\frac{5}{4}\omega(G)\rceil$.
\end{lemma}
\pf It holds clearly if $G$ is an induced subgraph of $F_{9}$.
We suppose that $\delta(G)\geq\lceil\frac{5}{4}\omega(G)\rceil$ by Lemmas~\ref{KM}.
If $Q_{x_1}=\emptyset$ or $Q_{x_{10}}=\emptyset$, then $G$ is a clique blowup of $F_{8}$ and we are done by Lemma~\ref{F_8}.
Suppose so that $Q_{x_1}\ne\emptyset$ and $Q_{x_{10}}\ne\emptyset$.
By Lemma~\ref{delta}, $|Q_{x_i}|\geq2$ for each $i\in\{2,7,8,9\}$.

If $Q_{x_4}=\emptyset$, then $Q_{5}=\emptyset$ because $G$ is connected and has no
clique cutsets. Now, $G$ is a clique blowup of $F_5$, and we are
done by Lemma~\ref{Petersen}. So, we further suppose, by symmetry,
that $Q_{x_4}\ne\emptyset$ and $Q_{x_5}\ne\emptyset$.
By Lemma~\ref{delta}, we have $|Q_{x_j}|\geq2$ for each $j\in\{3,4,5,6\}$.
Now, we can deduce that $|Q_k|\leq\omega(G)-2$ for all $k$.
Therefore, $\min\{|Q_v|, |Q_u|\}\ge 2$ for each edge $uv\in E(F_{10})$ with $|Q_{u, v}|=\omega(G)$.
Lemma~\ref{F_9} follows from Lemmas~\ref{2-blowup} and \ref{2-blowup'}.  \qed

\begin{lemma}\label{F_{10}'}
Let $G$ be a clique blowup of $F_{10}$. Then $\chi(G)\leq\lceil\frac{5}{4}\omega(G)\rceil$.
\end{lemma}
\pf It holds trivially if $G$ is an induced subgraph of $F_{10}$. We suppose $\delta(G)\geq\lceil\frac{5}{4}\omega(G)\rceil$ by Lemma~\ref{KM}.

If $Q_{x_4}=\emptyset$, then $Q_{x_5}=\emptyset$ as otherwise $Q_{x_6}$ is a clique cutset
of $G$ or $G$ is a complete graph. Now, $G$ is a clique blowup of $F_5$, and the lemma holds
by Lemma~\ref{Petersen}. So, suppose by symmetry that $Q_{x_4}\ne\emptyset$ and
$Q_{x_2}\ne\emptyset$. With the same argument, we may suppose that $Q_{x_5}\ne\emptyset$ and $Q_{x_1}\ne\emptyset$. Since $\delta(G)\geq\lceil\frac{5}{4}\omega(G)\rceil\geq\omega(G)+1$,
by Lemma~\ref{delta}, we have that $|Q_{x_i}|\geq2$ for each $i\in\{1,2,\cdots,7\}$. Then, $|Q_{x_k}|\leq\omega(G)-2$ for all $k$, and so $\min\{|Q_v|, |Q_u|\}\geq2$ for each $uv\in E(F_9)$ with $|Q_{uv}|=\omega(G)$.
Lemma~\ref{F_{10}'} follows from Lemmas~\ref{2-blowup} and \ref{2-blowup'}. \qed

\begin{lemma}\label{F_{11}}
Let $G$ be a clique blowup of $F_{11}$. Then $\chi(G)\leq\lceil\frac{5}{4}\omega(G)\rceil$.
\end{lemma}
\pf It holds clearly if $G$ is an induced subgraph of $F_{11}$.
We suppose that $\delta(G)\geq\lceil\frac{5}{4}\omega(G)\rceil$
by Lemma~\ref{KM}.

If $Q_{x_{11}}=\emptyset$, then $G$ is a clique blowup of $F_{8}$. If $Q_{x_8}=\emptyset$ or $Q_{x_9}=\emptyset$,
then $G$ is a clique blowup of $F_{10}$. If $Q_{x_4}=\emptyset$, then $Q_{x_5}=\emptyset$ since $G$ is connected, imperfect and has no clique cutsets, and $G$ is a clique blowup of $F_5$.
Therefore, by Lemmas~\ref{Petersen}, \ref{F_8} and \ref{F_{10}'}, we suppose
by symmetry  that $Q_{x_i}\ne\emptyset$ for each $i\in\{4,5,8,9,11\}$.

By Lemma~\ref{delta}, $|Q_{x_j}|\geq2$ for each $j\in\{2,3,\cdots,7\}$.
Now $|Q_k|\leq\omega(G)-2$ for all $k$, and so $\min\{|Q_v|,|Q_u|\}\geq2$ for each edge $uv\in E(F_{11})$ with $|Q_{u, v}|=\omega(G)$.
Lemma~\ref{F_{11}} follows from Lemmas~\ref{2-blowup} and \ref{2-blowup'}. \qed

\begin{lemma}\label{F_{12}}
Let $G$ be a clique blowup of $F_{12}$. Then $\chi(G)\leq\lceil\frac{5}{4}\omega(G)\rceil$.
\end{lemma}
\pf It holds if $G$ is an induced subgraph of $F_{12}$.
We still suppose that $\delta(G)\geq\lceil\frac{5}{4}\omega(G)\rceil$
by Lemma~\ref{KM}.

If $Q_{x_{1}}=\emptyset$ or $Q_{x_{10}}=\emptyset$ then $G$ is a clique blowup of $F_{11}$. If $Q_{x_{11}}=\emptyset$ then $G$ is a clique blowup of $F_{9}$. If $Q_{x_4}=\emptyset$,
then $Q_{x_5}=\emptyset$ since $G$ is connected, imperfect and has no clique cutsets, and so
$G$ is a clique blowup of $F_5$.
So, we suppose, by Lemmas~\ref{Petersen}, \ref{F_9} and \ref{F_{11}}, that
$Q_{x_{i}}\ne\emptyset$ for each $i\in\{1,4,5,10,11\}$.

By Lemma~\ref{delta}, $|Q_{x_j}|\geq2$ for each $j\in\{2,3,\cdots,9\}$.
Now, $|Q_{x_k}|\leq\omega(G)-2$ for all $k$, and $\min\{|Q_v|, |Q_u|\}\geq2$ for each edge $uv\in E(F_{12})$ with $|Q_{uv}|=\omega(G)$. Lemma~\ref{F_{12}} follows from Lemmas~\ref{2-blowup} and \ref{2-blowup'}. \qed

\medskip

Up to now, we have proved the first conclusion of Theorem~\ref{bull-3}, i.e., $\chi(G)\leq\lceil\frac{5}{4}\omega(G)\rceil$ for all $(C_4, F, H)$-free graphs. In the following Lemmas~\ref{C_5}$\sim$\ref{C_7'}, we show that we can do better by excluding some further configurations.

\begin{lemma}\label{C_5}
Let $G$ be a $C_5^2$-free clique blowup of $F_2$ or $F_5$.
Then $\chi(G)\leq\omega(G)+1$.
\end{lemma}
\pf It holds clearly if $G$ is an induced subgraph of $F_{2}$ or $F_5$. Let $G$ be a clique blowup of $F_2$ or $F_5$,
and suppose that $\delta(G)\geq\omega(G)+1$.

Let $S$ be a stable set of $G$. If $S$ is a good stable set,
then $\chi(G)\leq\chi(G-S)+1\leq\omega(G)+1$ by induction.
If $G-S$ is perfect, then $\chi(G)\leq\chi(G-S)+1\leq\omega(G)+1$.
We say a graph $H$ is {\em good} if $H$ has a stable set $S$ such that
either $S$ is a good stable set or $H-S$ is perfect.
Next, we show that $G$ is good, and thus satisfies $\chi(G)\leq\omega(G)+1$ by induction.
Suppose that
\begin{equation}\label{G-S-imperfect}
\mbox{$G-S$ is imperfect for any stable set $S$.}
\end{equation}

\begin{claim}\label{F_2+}
If $G$ is a clique blowup of $F_2$, then $G$ is good.
\end{claim}
\pf It is certain that $|Q_{x_i}|\geq2$ for each $i\in\{1,4,5\}$, since $G-\{x_i\}$ cannot be perfect. Since $G$ is $C_5^2$-free, we have that $\min\{|Q_{x_2}|, |Q_{x_3}|\}\le 1$,
and $\min\{|Q_{y_j}|, |Q_{z_j}|\}\le 1$ for each $j$.
Without loss of generality, we suppose that $|Q_{x_2}|\geq|Q_{x_3}|$.
For $j\in\{1,2,\cdots,t\}$, let $w_j\in\{y_j, z_j\}$ such that
$|Q_{w_j}|=\max\{|Q_{y_j}|,|Q_{z_j}|\}$. Then,
$Q_{x_2}\cup\{x_5, w_1, w_2,\cdots, w_t\}$ is a good stable set of $G$.
This proves Claim~\ref{F_2+}. \qed

Next, we suppose that $G$ is a clique blowup of $F_5$.
Note that $F_3$ and $F_4$ are both induced subgraphs of $F_5$.
We first discuss the cases that $G$ is a clique blowup of $F_3$ or $F_4$.

\begin{claim}\label{F_3+}
Let $G$ be a clique blowup of $F_3$ or $F_4$. Then, $G$ is good.
\end{claim}
\pf First suppose that $G$ is a clique blowup of $F_3$. Following from (\ref{imperfect000}), we have that $Q_{x_i}\ne\emptyset$ for all $i$
as otherwise $G$ is a clique blowup of $F_2$ and we are done.
By Lemma~\ref{delta}, $|Q_{x_j}|\geq2$ for each $j\in\{1,3,4,6\}$.
Since $G$ is $C_5^2$-free, we have that $\min\{|Q_{x_2}|, |Q_{x_5}|\}=\min\{|Q_{x_2}|,|Q_{x_7}|\}=\min\{|Q_{x_5}|,|Q_{x_{10}}|\}=\min\{|Q_{x_7}|, |Q_{x_{10}}|\}=1$.
Then, either $Q_{x_2, x_{10}}=\{x_2, x_{10}\}$ or $Q_{x_5, x_7}=\{x_5, x_7\}$, and so $G$ is good because both $G-\{x_2, x_{10}\}$ and $G-\{x_5, x_7\}$ are
perfect, which contradicts (\ref{G-S-imperfect}).

Suppose now that $G$ is a clique blowup of $F_4$. Then, $Q_{x_i}\ne\emptyset$ for all $i$
as otherwise $G$ is a clique blowup of $F_3$ and we are done.
By Lemma~\ref{delta}, $|Q_{x_j}|\geq2$ for each $j\in\{2,3,4,5,8,9\}$.
Since $G$ is  $C_5^2$-free, we have that $Q_{x_k}=\{x_k\}$ for $k\in\{1, 7, 10\}$. Now,
$\{x_1,x_7,x_{10}\}$ is a stable set and $G-\{x_1,x_7,x_{10}\}$ is perfect, which contradicts (\ref{G-S-imperfect}).
This proves Claim~\ref{F_3+}. \qed

\medskip

We now suppose that $G$ is a nonempty clique blowup of $F_5$ by Claim~\ref{F_3+}. By (\ref{G-S-imperfect}), we have that
\begin{equation}\label{not-stable}
\mbox{$Q_{u, v, w}$ is not a stable set}
\end{equation}
for any three distinct vertices $u,v,w\in V(F_5)$. Particularly, $N_G(Q_{x_i})$ is not a stable set, and so $|Q_{x_i}|\leq\omega(G)-2$ for all $i$.

Since $G$ is $C_5^2$-free, we suppose, without loss of generality, that $Q_{x_1}=\{x_1\}$.
Since none of $G[Q_{x_2, x_3, x_4, x_9, x_{10}}]$, $G[Q_{x_3, x_4, x_5, x_7, x_{8}}]$, and
$G[Q_{x_6, x_7, x_8, x_9, x_{10}}]$ can be a $C_5^2$, and since $N_G(\{x_1\})$ is not a stable set by (\ref{not-stable}),
we may suppose by symmetry that $Q_{x_{3}}=\{x_3\}$.
Now, we have that $|Q_{x_i}|\geq2$ for $i\in\{8, 9, 10\}$ by (\ref{not-stable}), and so $\min\{|Q_{x_6}|, |Q_{x_7}|\}=1$ since $G[Q_{x_6, x_7, x_8, x_9, x_{10}}]$ is not a $C_5^2$.

If $|Q_{x_7}|=1$, then $|Q_{x_4}|\geq2$ to avoid $Q_{x_1, x_4, x_7}$ being stable, and
$Q_{x_5}=\{x_5\}$ to avoid a $C_5^2$ on $Q_{x_4, x_5, x_8, x_9, x_{10}}$. But then,
$\{x_2, x_8, x_{10}\}$ is a good stable set of $G$ as $\omega(G[N_G[\{x_1, x_3, x_5, x_7\}]])<\omega(G)$.

Suppose that $|Q_{x_6}|=1$.Then, $|Q_{x_5}|\geq2$ for avoiding a stable set on $Q_{x_3, x_5, x_6}$, and $Q_{x_4}=\{x_4\}$ for avoiding a $C_5^2$ on $Q_{x_4, x_5, x_8, x_9, x_{10}}$. But then,
$\{x_2, x_8, x_{10}\}$ is a good stable set of $G$ as $\omega(G[N_G[\{x_1,x_3,x_4,x_6\}]])<\omega(G)$. This completes the proof of Lemma~\ref{C_5}. \qed

\begin{lemma}\label{clique}
Let $G$ be a clique blowup of $F_{12}$ such that $G$ has no good stable sets, and $Q_{x_8, x_9, x_{10}, x_{11}}$ meets no maximum clique of $G$. Then, $G[Q_{x_1, x_2, \cdots, x_7}]$ is isomorphic to $C_7^{\omega(G)/2}$. And, if $|Q_{x_i}|\leq1$ for each $i\in\{8, 9, 10, 11\}$ and every proper induced subgraph $G'$ satisfies $\chi(G')\leq\lceil 7\omega(G')/6\rceil$, then $\chi(G)\leq\lceil 7\omega(G)/6\rceil$.
\end{lemma}
\pf Since $Q_{x_8, x_9, x_{10}, x_{11}}$ meets no maximum clique of $G$, we have that, for each $i\in\{1,2,\cdots,7\}$, $|Q_{x_i, x_{i+1}}|=\omega(G)$ for avoiding a good stable set $\{x_{i+2}, x_{i+4}, x_{i+6}\}$, where the subscript are taken modulo 7. Therefore, $G[Q_{x_1, x_2, \cdots, x_7}]$ is isomorphic to $C_7^{\omega(G)/2}$.

Let $t=\omega(G)/2$, and suppose that  $|Q_{x_i}|\leq1$ for each $i\in\{8, 9, 10, 11\}$.  If $t=1$ or 2, then, $\chi(G)\leq\lceil 5\omega(G)/4\rceil=\lceil7\omega(G)/6\rceil$ by Lemma~\ref{F_{12}}. So, suppose that $t\ge 3$.

If $t=3$, then we can construct a 7-coloring $\phi$ of $G$ as follows: $\phi(Q_{x_1})\subseteq \{1, 2, 3\}$, $\phi(Q_{x_2})\subseteq \{4, 5, 6\}$, $\phi(Q_{x_3})\subseteq \{1, 2, 7\}$, $\phi(Q_{x_4})\subseteq \{3, 4, 5\}$, $\phi(Q_{x_5})\subseteq \{1, 6, 7\}$, $\phi(Q_{x_6})\subseteq \{2, 3, 4\}$, $\phi(Q_{x_7})\subseteq \{5, 6, 7\}$, $\phi(Q_{x_8})\subseteq \{7\}$, $\phi(Q_{x_9})\subseteq \{4\}$, $\phi(Q_{x_{10}})\subseteq \{1\}$, and $\phi(Q_{x_{11}})\subseteq \{5\}$. So, $\chi(G)\leq7=\lceil7\omega(G)/6\rceil$.

Let $T$ be a subset of $V(G)$ obtained by taking $\min\{3,|Q_v|\}$ vertices from $Q_v$ for each vertex $v\in V(F_{12})$. Clearly, $\chi(G[T])\leq7$ as shown above.
Then, $\chi(G)\leq\chi(G-T)+\chi(G[T])\leq\lceil7\omega(G-T)/6\rceil+7\leq\lceil7(\omega(G)-6)/6\rceil+7
=\lceil7\omega(G)/6\rceil$. This proves Lemma~\ref{clique}. \qed

\begin{lemma}\label{C_5'}
Let $G$ be a connected $C_5^2$-free clique blowup of $F_{12}$ such that $\delta(G)\geq\omega(G)+1$ and $G$ has no good stable sets. If $G$ is not a clique blowup of $F_3$, then $G$ has a stable set $S$ such that $G-S$ is perfect, or $G[Q_{x_1, x_2, \cdots, x_7}]$ is isomorphic to $C_7^t$ for some  $t$  and $|Q_{x_i}|\leq1$ for each $i\in\{8,9,10,11\}$.
\end{lemma}
\pf Suppose that $G$ is not a clique blowup of $F_3$. Clearly, $G$ is not a clique blowup of $F'_2$.  We divide the proof into several claims.
\begin{claim}\label{F_8-}
	If $G$ is a clique blowup of $F_{8}$, then the lemma holds.
\end{claim}
\pf For $i\in\{1,2,7\}$, since $F_8-x_i$ is perfect, we may suppose $|Q_{x_i}|\geq2$ as otherwise $G-\{x_i\}$ is perfect and so we are done. If $|Q_{x_8}|\geq 2$, then  $|Q_{x_6}|\geq2$ by Lemma~\ref{delta}, and so $G[Q_{x_1 x_2, x_6, x_7, x_8}]$ is a $C_5^2$. This shows that $|Q_{x_8}|\leq1$.  Similarly, $|Q_{x_9}|\leq1$. Since $G$ has no clique cutsets by (\ref{imperfect000}), we may suppose that $Q_{x_{4}}\ne\emptyset$ and $Q_{x_5}\ne\emptyset$ as otherwise $G$ is a clique blowup of $F'_2$. Again by Lemma~\ref{delta}, we have  $|Q_{x_j}|\geq2$ for $j\in\{3,4,5,6\}$, and so $|Q_{x_i}|\geq2$ for $i\in\{1, 2, \cdots, 7\}$. Recall that  $\max\{|Q_{x_8}|, |Q_{x_9}|\}\le 1$. Now, $Q_{x_8, x_9}$ meets no maximum clique of $G$, and  by Lemma~\ref{clique}, $G[Q_{x_1, \cdots, x_7}]$ is a $C_7^{\omega(G)/2}$. This proves Claim~\ref{F_8-}. \qed

\begin{claim}\label{F_9-}
	If $G$ is a clique blowup of $F_{9}$, then the lemma holds.
\end{claim}
\pf If $Q_{x_i}=\emptyset$ for some $i\in\{1,8,9,10\}$ then $G$ is a clique blowup of $F_8$ since $G$ is a connected imperfect graph without clique cutsets by (\ref{imperfect000}).  If $Q_{x_j}=\emptyset$ for some $j\in\{2,3,\cdots,7\}$, then $G$ is a clique blowup of $F_3$. So, we suppose that $Q_{x_i}\ne \emptyset$ for all $i$, and so $|Q_{x_k}|\geq2$ for each $k\in\{2,3,\cdots,9\}$ by Lemma~\ref{delta}. Then, $|Q_{x_1}|=|Q_{x_{10}}|=1$ since $G$ is $C_5^2$-free, and so $\{x_1,x_{10}\}$ is a stable set with $G-\{x_1,x_{10}\}$ perfect. This proves Claim~\ref{F_9-}.  \qed

\begin{claim}\label{F_{10}-}
	If $G$ is a clique blowup of $F_{10}$, then the lemma holds.
\end{claim}
\pf By Claim~\ref{F_8-}, we may suppose that $Q_{x_9}\ne\emptyset$ and $Q_{x_{11}}\ne\emptyset$ as otherwise $G$ is a clique blowup of $F_{8}$. We may also suppose $Q_{x_i}\ne\emptyset$ for each $i\in\{1,2,4,5\}$ as otherwise $G$ is a clique blowup of $F'_2$.
Then, $|Q_{x_j}|\geq2$ for each $j\in\{1,2,\cdots,7\}$ by Lemma~\ref{delta}, and $|Q_{x_9}|=|Q_{x_{11}}|=1$ since $G$ is $C_5^2$-free, and so  $Q_{x_9, x_{11}}$ meets no maximum clique of $G$. Now, $G[Q_{x_1, \cdots, x_7}]$ is $C_7^{\omega(G)/2}$ by Lemma~\ref{clique}.
This proves Claim~\ref{F_{10}-}. \qed

\begin{claim}\label{F_{11}-}
If $G$ is a clique blowup of $F_{11}$, then the lemma holds.
\end{claim}
\pf By (\ref{imperfect000}), we may suppose $Q_{x_i}\ne\emptyset$ for each $i\in\{1,2,\cdots,9,11\}$ as otherwise $G$ is a clique blowup of $F_3$, or $F_8$, or $F_{10}$, and we are done by Claims~\ref{F_8-} and \ref{F_{10}-}. Then, $|Q_{x_j}|\geq2$ for each $j\in\{2\,\cdots,7\}$ by Lemma~\ref{delta}, and $|Q_{x_{11}}|=1$ since $G$ is $C_5^2$-free. If $|Q_{x_1}|=1$, then $\{x_1,x_{11}\}$ is a stable set with  $G-\{x_1,x_{11}\}$ perfect. If $|Q_{x_1}|\geq2$, then $\max\{|Q_{x_8}|,|Q_{x_9}|\}\leq1$ since $G$ is $C_5^2$-free. Hence, $Q_{x_8,x_9,x_{11}}$ meets no maximum clique of $G$, and now $G[Q_{x_1, \cdots, x_7}]$ is $C_7^{\omega(G)/2}$ by Lemma~\ref{clique}. This proves Claim~\ref{F_{11}-}.\qed

\medskip

Now, we may suppose that $G$ is a nonempty clique blowup of $F_{12}$, as otherwise the statement holds by Claims~\ref{F_8-}$\sim$\ref{F_{11}-}.
Then, by Lemma~\ref{delta},  $|Q_{x_j}|\geq2$ for each $j\in\{2,3,4,5,6,7,8,9\}$, and so $|Q_{x_1}|=|Q_{x_{10}}|=|Q_{x_{11}}|=1$ since $G$ is $C_5^2$-free. Now, $\{x_1,x_{10},x_{11}\}$ is a stable set and $G-\{x_1, x_{10}, x_{11}\}$ is perfect. This completes the proof of Lemma~\ref{C_5'}. \qed


\begin{lemma}\label{C_5''}
Let $G$ be a $C_5^2$-free clique blowup of $F_{12}$. Then, $\chi(G)\leq\lceil\frac{7}{6}\omega(G)\rceil$.
\end{lemma}
\pf It holds trivially if $G$ is an induced subgraph of $F_{12}$. We may suppose by induction that $\delta(G)\geq\lceil\frac{7}{6}\omega(G)\rceil$, $G$ has no good stable sets, and has no stable set $S$ such that $G-S$ is perfect.

By Lemma~\ref{C_5'},  either $G$ is a clique blowup of $F_3$, or $G[Q_{x_1, x_2, \cdots, x_7}]$ is isomorphic to $C_7^t$ for some $t$ and $|Q_{x_i}|\leq1$ for $i\in\{8,9,10,11\}$. Recall that $F_3$ is an induced subgraph of $F_5$. Then, Lemma~\ref{C_5''} follows from Lemmas~\ref{C_5} and \ref{clique}. \qed

\begin{lemma}\label{C_7'}
Let $G$ be a $(C_5^2,C_7^4)$-free clique blowup of $F_{12}$. Then $\chi(G)\leq\omega(G)+1$.
\end{lemma}
\pf The statement holds easily if $G$ is an induced subgraph of $F_{12}$. Suppose  by induction that $\delta(G)\geq\omega(G)+1$, $G$ has no good stable sets, and has no stable set $S$ such that $G-S$ is perfect.

By Lemma~\ref{C_5'},  either $G$ is a clique blowup of $F_3$, or $G[Q_{x_1, x_2, \cdots, x_7}]$ is isomorphic to $C_7^{\omega(G)/2}$ and $|Q_{x_i}|\leq1$ for each $i\in\{8,9,10,11\}$. If $G$ is a clique blowup of $F_3$, we are done by Lemma~\ref{C_5}. So, suppose that $G$ satisfies the latter.
Since $G$ is $C_7^4$-free, we have that $\omega(G)\leq6$, and so    $\chi(G)\leq\lceil\frac{7}{6}\omega(G)\rceil=\omega(G)+1$ by Lemma~\ref{C_5''}. This completes the proof of Lemma~\ref{C_7'}. \qed

\medskip

Now, we are ready to prove Theorem~\ref{bull-3}.

\noindent\textbf{{\em Proof of Theorem~\ref{bull-3}}:}  Let $H$ be a $P_7$ or a fork$^+$, and let $G$ be a ($C_4, H$, bull)-free graph. We may suppose that $G$ is connected, imperfect and has no clique cutsets or universal cliques. By Theorems~\ref{P_7}, \ref{fork^+}, and Lemma~\ref{perfect'}, $G$ is a nonempty clique blowup of an  imperfect graph, where the graph is in $\mathfrak{F}$ if $H$ is a $P_7$, and is in $\mathfrak{F'}$ if $H$ is a fork$^+$.

Lemmas~\ref{Petersen}, \ref{odd-hole'}, \ref{F_2} and \ref{F_{12}},
assert that $\chi(G)\leq\lceil\frac{5}{4}\omega(G)\rceil$, and Lemmas~\ref{odd-hole'},
\ref{C_5} and \ref{C_5''} assert that $\chi(G)\leq\lceil\frac{7}{6}\omega(G)\rceil$ if $G$ is ($C_4, H$, bull, $C_5^2$)-free.

By Lemma~\ref{C_5} and \ref{C_7'}, we have that $\chi(G)\leq\omega(G)+1$ if $G$ is ($C_4, P_7$, bull, $C_5^2, C_7^4$)-free.

Suppose finally that $G$ is ($C_4$, fork$^+$, bull, $C_5^2, C_7^4$)-free. Then, $G$ is a nonempty clique blowup of a graph $F\in \{F'_2, F_3, F_4, F_5\}\cup\{C_k\;|\; k\ge 5\}$. If $F\in \{F'_2, F_3, F_4, F_5, C_5\}$ then  $\chi(G)\leq\omega(G)+1\le \lceil\frac{9}{8}\omega(G)\rceil$ by Lemma~\ref{C_5}. If $F=C_7$ then we are done by Lemma~\ref{C_7'}. If $F=C_k$ for some $k=6$ or $k\ge 8$, then by Lemmas~\ref{perfect'} and \ref{odd-hole'}, $\chi(G)\leq\max_{q\ge 4}\lceil\frac{2q+1}{q}\omega(G)\rceil\le  \lceil\frac{9}{8}\omega(G)\rceil$. This completes the proof of
of Theorem~\ref{bull-3}. \qed

\medskip

\noindent{\bf Remark:} Let $\Delta(G)$ denote the maximum degree of $G$. Brooks \cite{B1941} showed that if $G$ is a graph with $\Delta(G)\geq3$ and $\omega(G)\leq\Delta(G)$, then $\chi(G)\leq\Delta(G)$. Reed \cite{R1998} conjectured that every graph $G$ satisfies  $\chi(G)\leq\lceil\frac{\omega(G)+\Delta(G)+1}{2}\rceil$. This conjecture is still widely open. Gernet {\em et al} \cite{GR2008} showed that Reed's conjecture holds for graphs $G$ with $\chi(G)\leq\omega(G)+2$, and Karthick {\em et al} \cite{KM2018} showed that Reed's conjecture holds for graphs $G$ with $\chi(G)\leq\lceil\frac{5}{4}\omega(G)\rceil$. Combining these two conclusions with Corollary~\ref{hammer-3} and Theorem~\ref{bull-3},  we have immediately the following theorem.

\begin{theorem}\label{Reed's Conj-2}
Let $F$ be a $P_7$ or a fork$^+$, and let $H$ be a bull or a hammer. Then, Reed's conjecture holds for  $(C_4, F, H)$-free graphs.
\end{theorem}

\end{document}